\documentclass{article}
\usepackage{amsmath,amsthm,amsfonts}
\newtheorem{definition}{Definition}
\newtheorem{theorem}{Theorem}
\newtheorem{proposition}{Proposition}
\theoremstyle{remark}
\newtheorem{remark}{Remark}

\begin{document}
\title{Foundations of Real Analysis and Computability Theory in Non-Aristotelian Finitary Logic}
\author{Radhakrishnan Srinivasan\thanks{ISL-Research, IBM India Software Labs, Ozone-1 (Ground Floor), SP Infocity,
Pune-Saswad Road, Fursungi, Pune 412308. E-Mail: sradhakr@in.ibm.com} and
H.~P.~Raghunandan\thanks{IBM Center for Advanced Studies, IBM India Software Labs, Ozone-1 (Ground Floor),
SP Infocity, Pune-Saswad Road, Fursungi, Pune 412308. E-Mail: hpraghu@in.ibm.com}}
\date{}
\maketitle

\begin{abstract}
This paper outlines new paradigms for real analysis and computability theory in the recently proposed
non-Aristotelian finitary logic~(NAFL). Constructive real analysis in NAFL~(NRA) is accomplished by a translation
of diagrammatic concepts from Euclidean geometry into an extension~(NPAR) of the NAFL version of Peano Arithmetic~(NPA).
Such a translation is possible because NPA proves the existence of every infinite proper class of
natural numbers that is definable in the language of NPA. Infinite sets are not permitted in NPAR
and quantification over proper classes is banned; hence Cantor's diagonal argument cannot be
legally formulated in NRA, and there is no `cardinality' for any collection (`super-class') of real numbers. Many of the
useful aspects of classical real analysis, such as, the calculus of Newton and Leibniz, are justifiable in NRA.
But the paradoxes, such as, Zeno's paradoxes of motion and the Banach-Tarski paradox, are resolved because NRA
admits only closed super-classes of real numbers; in particular, open/semi-open intervals of real numbers are not
permitted. The NAFL version of computability theory~(NCT) rejects Turing's argument for the undecidability of the
halting problem and permits hypercomputation. Important potential applications of NCT are in the areas of quantum
and autonomic computing.
\end{abstract}

\section{Introduction to NAFL}\label{intro}
The basic description of non-Aristotelian finitary logic~(NAFL) given in \cite{ijqi,1166,bcp,acs} is
outlined below for the sake of completeness. The language, well-formed formulae and rules of inference of
NAFL theories are formulated in exactly the same manner as in classical
first-order predicate logic with equality~(FOPL), where we shall assume, for convenience,
that natural deduction is used; however, there are key differences and restrictions
imposed by the requirements of the Main Postulate of NAFL, which is defined below.
In NAFL, truths for \emph{formal propositions} can exist \emph{only} with respect
to axiomatic theories. There are no absolute truths in just the \emph{language} of
an NAFL theory, unlike classical/intuitionistic/constructive logics. There do exist
absolute (metamathematical, Platonic) truths in NAFL, but these are truths
\emph{about} axiomatic theories and their models. As in FOPL, an NAFL theory T
is defined to be consistent if and only if T has a model, and a proposition
$P$ is undecidable in T if and only if neither $P$ nor its negation
$\neg P$ is provable in T. The truth definition given below in the Main Postulate is the
heart of NAFL and its vital importance is underscored by Proposition~\ref{p1}, which follows it.
The existing mindset that ``truth definitions are not part of logic'' most emphatically does not
hold of NAFL.

\subsection{The Main Postulate of NAFL}\label{mp}
Let $P$ be a legitimate proposition of a consistent NAFL theory T. We will see in Sect.~\ref{tsps}
that legitimacy of a proposition requires it to be in the `theory syntax' of T.
If $P$ is provable/refutable in T, then $P$ is true/false with respect to T (henceforth
abbreviated as `true/false in T'); i.e., a model for T will assign $P$ to be true/false.
If $P$ is undecidable in a consistent NAFL theory T, then the Main Postulate provides the appropriate
truth definition as follows: $P$ is true/false in T if and only if $P$ is provable/refutable in an
\emph{interpretation} T* of T. Here T* is an axiomatic NAFL theory that, like T, temporarily resides
in the human mind and acts as a `truth-maker' for (a model of) T. The theorems of T* are precisely
those propositions that are assigned `true' in the NAFL model of T, which, unlike its classical
counterpart, is not `pre-existing' and is instantaneously generated by T*.
Obviously, the theorems of T must necessarily be included in those of T*.
Note that for a given consistent theory T, T* could vary in time according to the
free will of the human mind that interprets T; for example, T* could be T+$P$ or
T+$\neg P$ or just T itself at different times for a given human mind, or in the context
of quantum mechanics, for a given 'observer'. Further, T* could vary from one observer
to another at any given time; each observer determines T* by his or her own free will.
The essence of the Main Postulate is that $P$ is true/false in T if and only if it has been
\emph{axiomatically declared} as true/false by virtue of its provability/refutability
in T*. In the absence of any such axiomatic declarations, i.e., if $P$ is
undecidable in T* (e.g.~take T*=T), then $P$ is `neither true nor false' in T and
Proposition~\ref{p1} shows that consistency of T requires the laws of the excluded middle
and non-contradiction to fail in a non-classical model for T in which the superposition
$P \& \neg P$ is the case.

\begin{proposition}\label{p1}
Let $P$ be undecidable in a consistent \emph{NAFL} theory \emph{T}. Then $P \vee \neg P$
and $\neg (P \& \neg P)$ are not theorems of \emph{T}. There must exist a non-classical
model $\mathcal{M}$ for \emph{T} in which $P \& \neg P$ is the case.
\end{proposition}
For a proof of Proposition~\ref{p1}, see \cite{ijqi} or \cite{acs}; this proof, also reproduced in
Appendix~A, seriously questions the logical/philosophical basis for the law
of non-contradiction in both classical and intuitionistic logics.
The non-classical model $\mathcal{M}$ of Proposition~\ref{p1} is a superposition of two or
more classical models for T, in at least one of which $P$ is true and $\neg P$ in another. Here
`classical' or `non-classical' is used strictly with respect to the status of $P$, and `superposition'
means that \emph{all} the truths in each of the superposed classical models must hold in the non-classical
model $\mathcal{M}$. Thus if $P$ is true in one classical model and $\neg P$ in another, the superposed
state $P \& \neg P$ will hold when these two classical models are superposed to form $\mathcal{M}$.  In
$\mathcal{M}$, `$P$'~(`$\neg P$') denotes that `$\neg P$'~(`$P$') is not provable in T*, or
in other words, $\mathcal{M}$ expresses that neither $P$ nor $\neg P$ has been
axiomatically declared as (classically) true with respect to T; thus $P$, $\neg P$,
and hence $P \& \neg P$, are indeed (non-classically) true \emph{in our world}, according
to their interpretation in $\mathcal{M}$. Note also that $P$ and $\neg P$ are
\emph{classically} `neither true nor false' in $\mathcal{M}$, where `true' and `false' have
the meanings given in the Main Postulate. The quantum superposition principle is justified
by identifying `axiomatic declarations' of truth/falsity of $P$ in T (via its
provability/refutability in T*) with `measurement' in the real world~\cite{ijqi,1166,bcp,acs}.

Proposition~\ref{p1} is a metatheorem, i.e., it is a theorem
\emph{about} axiomatic theories. The concepts in Proposition~\ref{p1}, namely, consistency,
undecidability~(provability) and the existence of a non-classical model for a theory
and hence, quantum superposition and entanglement, are strictly metamathematical
(i.e., pertaining to semantics or model theory) and not formalizable in the
syntax of NAFL theories. An NAFL theory~T is either consistent or inconsistent,
and a proposition~$P$ is either provable or refutable or undecidable in T, i.e.,
the law of the excluded middle applies to these metamathematical truths.
Note that the existence of a non-classical model does not make T inconsistent or
paraconsistent, because T does not \emph{prove} $P \& \neg P$. However, one could assert that
the model theory for T requires the framework of a paraconsistent logic, in which the
non-classical models can be analyzed. NAFL is the only logic that correctly embodies
the philosophy of formalism~\cite{1166}; NAFL truths for formal propositions are
axiomatic, mental constructs with strictly no Platonic world required.

\subsection{Theory syntax and proof syntax}\label{tsps}
An NAFL theory T requires two levels of syntax, namely the `theory syntax' and the
`proof syntax'. The theory syntax consists of precisely those propositions that are
legitimate, i.e., whose truth in T satisfies the Main Postulate; obviously,
the axioms and theorems of T are required to be in the theory syntax. Further, one can
only add as axioms to T those propositions that are in its theory syntax. In particular,
neither $P \& \neg P$ nor its negation $P \vee \neg P$ is in the theory syntax when
$P$ is undecidable in T; this will be clear from the explanation given below. The proof syntax,
however, is classical because NAFL has the same rules of inference as FOPL; thus $\neg (P \& \neg P)$
is a valid deduction in the proof syntax and may be used to prove theorems of T. For
example, if one is able to deduce $A \Rightarrow P \& \neg P$ in the proof syntax of T
where $P$ is undecidable in T and $A$ is in the theory syntax, then one has proved $\neg A$
in T despite the fact that $\neg (P \& \neg P)$ is not a theorem (in fact not even a
legitimate proposition) of T. This is justified as follows: $\neg (P \& \neg P)$ may be
needed to prove theorems of T, but it does not follow in NAFL that
the theorems of T imply $\neg (P \& \neg P)$ if $P$ is undecidable in T.
Let $A$ and $B$ be undecidable propositions in the theory syntax of T.
Then $A \Rightarrow B$ (equivalently, $\neg A \vee B$) is
in the theory syntax of T if and only if $A \Rightarrow B$ is \emph{not}
(classically) deducible in the proof syntax of T. It is easy to check that if
$A \Rightarrow B$ is deducible in the proof syntax of T, then its (illegal) presence in the
theory syntax would force it to be a theorem of T, which is not permitted by
the Main Postulate. For, in a non-classical model $\mathcal{M}$ for T in which both $A$
and $B$ are in the superposed state (i.e., both $A \& \neg A$ and $B \& \neg B$ hold),
$A \& \neg B$ must be non-classically true, but theoremhood of $A \Rightarrow B$ will prevent
the existence of $\mathcal{M}$. If one replaces $B$ by $A$ in this result, one obtains the previous conclusion that
$\neg (A \& \neg A)$ is not in the theory syntax. For example, take $\mbox{T}_0$ to be the
null set of axioms. Then nothing is provable in $\mbox{T}_0$, i.e., every legitimate
proposition $P$ of $\mbox{T}_0$ is undecidable in $\mbox{T}_0$. Hence $P \& \neg P$ must
be satisfiable in a non-classical model for $\mbox{T}_0$; this is obviously true, for $P \& \neg P$ contradicts
only $P \vee \neg P$, which is not in the theory syntax (and hence cannot be a theorem) of $\mbox{T}_0$.
In particular, the proposition $(A \& (A \Rightarrow B))\Rightarrow B$, which is deducible in the proof
syntax of $\mbox{T}_0$ (via the \emph{modus ponens} inference rule), is not in the theory syntax; however,
if $A \Rightarrow B$ is not deducible in the proof syntax  of $\mbox{T}_0$,
then it is in the theory syntax. Note also that $\neg \neg A \Leftrightarrow A$ is
not in the theory syntax of $\mbox{T}_0$; nevertheless, the `equivalence' between
$\neg \neg A$ and $A$ holds in the sense that one can be replaced by the other
in every model of $\mbox{T}_0$, and hence in all NAFL theories. Indeed, in a non-classical
model for $\mbox{T}_0$, this equivalence holds in a non-classical sense and must be
expressed by a different notation~\cite{1166}.

\section{Infinite Sets Do Not Exist in NAFL Theories}\label{infsets}
Consider set theory with classes and assume that the only mathematical objects permitted to belong to classes
are sets. A class is specified as an extension of a definite property in the language of set theory. Let $P(x)$
be one such property of sets, in which all variables (bound and unbound) are restricted to be set variables and in which
finitely many independently defined constant class symbols may appear. We denote the class $C$ associated with
$P(x)$ as $\forall x[x \in C \Leftrightarrow P(x)]$ or more concisely, $C = \{x:P(x)\}$. 
Certain classes may be sets; a proper class is a collection that is deemed to be not a set and therefore
cannot belong to any class. The universe of a model for a theory is a nonempty class $U$ over which
all the variables of the theory range. The theory of finite sets~F, with classes, may be obtained by
removing the axiom of infinity~(AI) from either Zermelo's set theory with classes~(Z)~\cite[Chaps.~1~and~5]{Mcr},
or from G\"odel-Bernays theory~(GB)~\cite[pp.~73--78]{Coh}. That every property expressible in the language
determines a class is a theorem scheme of GB~\cite[pp.~76--77]{Coh} and in particular, of F; we denote this
as the theorem of comprehension for classes. The natural numbers in F are defined in the usual set-theoretic
sense, with $0$ as the empty set $\emptyset$ and $n+1 = n \cup \{n\}$. The class $N$ of all natural numbers is
defined by
\[
N=\{x: x \; \mbox{is a natural number}\}.
\]
Note that we have admitted ``$x$ is a natural number'' as a valid property in the language of F.
Consider the class $D$, defined as:
\begin{equation}\label{eq2}
D = \{x: \forall n(n \in \mbox{N}
\Rightarrow x \notin \wp^{(n)}(\emptyset))\}.
\end{equation}
Here $\wp^{(n)}(\emptyset)$ is the power set operation performed $n$ times on the null set $\emptyset$; we
let $\wp^{(0)}(\emptyset) = \emptyset$. Note that $D$ exists provably in F by the theorem of comprehension.
Consider the following proposition:
\begin{equation}\label{eq3}
D = \emptyset.
\end{equation}
According to the standard interpretation, (\ref{eq3}) is undecidable in F, for the following reason. Define a
model $\mathfrak{B}$ for F in which only finite sets exist; the universe $B$ of $\mathfrak{B}$ is given by
\begin{equation}\label{eq4}
B = \{x: \exists n(n \in \mbox{N} \; \& \;
x \in \wp^{(n)}(\emptyset))\}.
\end{equation}
In the interpretation of $\mathfrak{B}$, $D = \emptyset$ holds and hence (\ref{eq3}) is true. But in a model
in which infinite sets do exist, (\ref{eq3}) is false.
\begin{theorem}\label{t1}
The proposition $D=\emptyset$ defined in (\ref{eq3}) is required to be provable in \emph{F} by the Main Postulate
of \emph{NAFL} and Proposition~\ref{p1}. Consequently, infinite sets do not exist in the \emph{NAFL} version of \emph{F}.
\end{theorem}
\begin{proof}
We will outline the proof given in \cite[Sect.~3]{635}. Uniqueness of $D$ is deducible from the axiom of extensionality
for classes~(AE), which implies:
\begin{equation}\label{eq5}
D = \emptyset \; \vee \; D \ne \emptyset.
\end{equation}
The key point is that AE \emph{requires} (\ref{eq5}) (equivalent to $\neg (D=\emptyset \; \& \; D \ne \emptyset)$)
to be a theorem of F. But by Proposition~\ref{p1}, (\ref{eq5}) cannot be a theorem of F, given
that F is consistent and that (\ref{eq3}) is undecidable in F. We are forced to the conclusion
that one of these two ``givens'' must be false in NAFL. Either (\ref{eq3}) is undecidable in F, which is therefore
inconsistent, or else (\ref{eq3}) must be decidable in F.

What we have demonstrated is that if F is to be a consistent theory of NAFL, then we need a principle
of proof for either $D=\emptyset$ or for $D \ne \emptyset$. We will argue the
case for $D=\emptyset$. The axioms of F tacitly presume that a definite class $U$
is specified as the universe, because these formulas contain the universal quantifiers $\forall$ and
$\exists$ which may be thought of as referring to the universe.
However, to define $U$ we need to consider the formulas of F as meaningful because these assert
the existence of sets that must necessarily belong to $U$. Thus both $U$ and all the sets of F are
impredicatively defined. However, finite sets may also be defined predicatively, at least in principle,
merely by listing all of their elements; but this is not possible for infinite sets.
In conclusion, infinite sets can only be defined impredicatively, i.e.,
by self-reference via the universal class (which must again be defined impredicatively by reference to
infinite sets). The appropriate principle of proof in NAFL for $D=\emptyset$
is that such an impredicative definition of sets must be banned, for it generates undecidable propositions
that contradict the Main Postulate of NAFL.
\end{proof}
\begin{remark}
Note that Theorem~\ref{t1} is a metatheorem and its proof required Proposition~\ref{p1}. Classically, it is
not known whether $D=\emptyset$ is either provable or refutable or undecidable in F; the undecidability is
merely \emph{assumed}. What Theorem~\ref{t1} \emph{predicts} is that by the classical rules of inference, the theory
$\mbox{F} + (D \ne \emptyset)$ must be inconsistent. In other words, Theorem~\ref{t1} proves that F proves (\ref{eq3})
(assuming consistency of F), but does not provide a specific proof in F of (\ref{eq3}).
\end{remark}
\begin{remark}
The proof of Theorem~\ref{t1} illustrates a very important requirement of NAFL. Whenever an NAFL theory proves that
a class or a set must exist uniquely (e.g.\ via AE), then Proposition~\ref{p1} requires that such uniqueness must be
enforced with respect to that \emph{theory}, rather than just with respect to \emph{models} of that theory. Thus the
theory is also required to provide a unique \emph{construction} for the class or set in question. E.g.\ NAFL requires
the theory F to specify the class $D$ uniquely; it is not valid in NAFL to assert that either $D = \emptyset$ or
$D \ne \emptyset$ can hold in models of F, for then uniqueness is lost with respect to F. One can immediately
conclude that non-Euclidean geometries are inconsistent in NAFL~\cite{635,1255}, for Euclid's first four postulates
prove that there must exist a \emph{unique} straight line between any two distinct points. The uniqueness here must
again be enforced with respect to a theory comprising the first four postulates, which implies that such a theory must
provide a \emph{construction} for the said straight line via provability of the fifth postulate.
\end{remark}

\section{Infinite proper classes must provably exist in NAFL theories with infinite domains}
\begin{proposition}\label{p2}
If \emph{T} is a consistent \emph{NAFL} theory that proves the existence of infinitely many objects identifiable
by a property $P(x)$ in the language of \emph{T} (e.g.\ $P(x)$ could be `$x=x$' with \emph{T}
taken as Peano Arithmetic) then \emph{T} must prove the existence of the infinite proper class $\{x: P(x)\}$.
\end{proposition}
Note that Proposition~\ref{p2} addresses only \emph{infinite} classes; there is no such corresponding result for
finite classes (or sets). In a sense the axioms for the existence of infinite classes in the theory F of
Sect.~\ref{infsets} are redundant; by Proposition~\ref{p2}, they must be provable in the NAFL version of F, and
the language of F must necessarily include the required properties that specify the infinite classes.
\begin{proof}
The simplest proof is as follows. Consider the NAFL version of Peano Arithmetic~(NPA; see Appendix~B), which is
equivalent to F. From Theorem~\ref{t1}, it follows in NAFL that nonstandard models of NPA cannot exist. For the
Main Postulate of NAFL (see Sect.~\ref{intro}) requires that models of NPA must necessarily be specified by its
interpretation NPA*, which is also an NAFL theory that denies the existence of infinite sets (by Theorem~\ref{t1}).
But without infinite sets, nonstandard models of arithmetic cannot be specified and therefore cannot exist in NAFL. It
immediately follows that the only model of NPA is the standard model and hence, by the completeness theorem
of first-order logic, NPA must prove the existence of the domain of standard natural numbers specified by the
class $N$ defined in Sect.~\ref{infsets} (in the language of NPA, $N$ may be specified as $\{x: x = x \}$). This
result is immediately generalized and Proposition~\ref{p2} follows.
\end{proof}
\begin{remark}\label{rmnon}
Note that the completeness theorem of classical first-order logic, which speaks \emph{about} axiomatic theories,
is taken for granted in NAFL as a metamathematical principle, but is not provable in any NAFL theory because of
the infinitary reasoning required; in fact the notion of `axiomatic theory' is not even formalizable in NAFL theories.
The nonexistence of nonstandard models of NPA also follows from the arguments given in \cite[Sect.~1]{1166},
which raise serious questions about the existence of nonstandard integers even within classical logic.
Indeed, undecidability of the G\"odel sentence for NPA immediately produces a contradiction similar to
that discussed in the proof of Theorem~\ref{t1}. It follows as a corollary that consistency of NPA
demands its \emph{completeness}, i.e., every sentence of NPA must be either provable or refutable in NPA.
This result, together with Theorem~\ref{t1}, means that one no longer has the freedom to invoke infinite
sets to prove purely arithmetical propositions, as is done classically. Obviously, G\"odel's incompleteness
theorems and Turing's argument for the undecidability of the halting problem fail in NAFL, as discussed in
\cite{1166} and \cite{acs}. Another consequence of the nonexistence of nonstandard models is that an
`arbitrary but fixed' constant does not exist in NAFL; it is the same as a free variable and assumed to be
universally quantified with respect to a \emph{standard} domain (e.g.\ natural numbers). The Main Postulate
requires that the free (or universally quantified) variable is in a superposed state of all possible
values~\cite{1166} and assumes a particular value when and only when it is specified by the human mind that
interprets the theory in question. 
\end{remark}
\begin{remark}\label{rmcl}
An intuitive and useful explanation for the validity of Proposition~\ref{p2} is as follows. In NAFL, the
process of counting one or finitely many natural numbers at a time will not exhaust `all' the natural numbers.
A simple induction tells us that such a counting process will always leave out infinitely many natural
numbers, and so will never be complete. Yet NPA obviously specifies `all' natural numbers; it can only do so by
accessing (at some unspecified point in the counting process) infinitely many natural numbers \emph{at the same time}.
Thus NPA must prove the existence of the infinite proper class of natural numbers, which may be thought of as the
\emph{simultaneous} existence of `all' natural numbers. A class is identified by all of its elements, as required
by AE. When a theory `constructs' all the elements of an infinite class, it must unavoidably construct the class itself.
Proposition~\ref{p2} is a metatheorem; it proves that T must prove the existence of the infinite class $\{x: P(x)\}$,
but does not provide the said proof in T. The proof in T that infinitely many elements exist that
satisfy $P(x)$ may also be taken as the proof in T of the existence of $\{x: P(x)\}$. For example,
if $P(x)$ expresses that `$x$ is a natural number' (say, via $P(x) \equiv `x = x'$) in the theory NPA,
one might take the proof in NPA of the existence of the infinite class $N = \{x: P(x)\}$ as that of the
proposition $\forall n \exists m (m > n)$.
\end{remark}
\begin{remark}\label{ae}
It is clear from Proposition~\ref{p2} and Remark~\ref{rmcl} that an NAFL theory T always identifies an infinite
class $\{x: P(x)\}$ by all of its elements. In fact $\{x: P(x)\}$ can \emph{only} be so identified as an infinite
class, i.e., in the extensional sense. Thus the axiom of extensionality for classes, in those instances
where infinite classes are involved,  is built into the very definition of the word `class' and must also be
provable in every NAFL theory T that proves the existence of infinitely many objects.
\end{remark}
\begin{proposition}\label{p3}
`Arbitrary' infinite classes and quantification over infinite classes (which are equivalent; see Remark~\ref{rmnon})
are not permitted in \emph{NAFL} theories. An infinite (proper) class must always be a constant that is specified
constructively in \emph{NAFL} theories, via a mapping to the class $N$ of all natural numbers (which is assumed to
be a constructive specification). This mapping must be shown to exist even when one is unable to specify a
general formula for it (e.g.\ the class of all prime natural numbers).
\end{proposition}
Proposition~\ref{p3} emphasizes the difference between an infinite set (which is a mathematical object
in classical logic) and an infinite proper class. Quantification over proper classes amounts to treating
these as sets.
\begin{proof}
From Proposition~\ref{p2} and Remark~\ref{rmnon}, it is clear that G\"odel's incompleteness theorems (and
Turing's argument for the unsolvability of the halting problem) must fail in NAFL. To justify such a failure,
one must impose a ban on quantification over proper classes. The proofs of G\"odel's theorems (as well as those of all
versions of Turing's argument for the undecidability of the halting problem~\cite{acs}) always
require quantification over infinitely many infinite proper classes. The said quantification in various versions
of the proofs of these theorems appears either directly (e.g.\ via Cantor's diagonalization argument) or indirectly,
e.g.\ when one encodes self-referential sentences into Peano Arithmetic using G\"odel numbering. In other words, G\"odel's
theorems and Turing's argument serve as \textit{reductio ad absurdum} proofs of the infinitary nature, and hence,
invalidity in NAFL, of quantification over proper classes. Indeed, the very fact that nonstandard models of arithmetic
follow as a logical consequence of G\"odel's theorems and Turing's argument is proof of their infinitary nature.
\end{proof}

\section{Foundations of Real Analysis in NAFL}\label{ra}
Extend NPA to the NAFL theory NPAR in which the integers~($Z$), rationals~($Q$) and sequences of
rationals~($\langle q_n \rangle$), including Cauchy sequences, are defined. See Appendix~B for the definitions of
NPA and NPAR. By Propositions~\ref{p2}~and~\ref{p3}, NPAR must prove the existence of every Cauchy sequence of rationals
definable by formulas (`properties') in the language of NPAR. The usual extension of NPAR to treat real numbers as
equivalence classes of such Cauchy sequences is not possible, for direct quantification over reals is banned.
Neither can we talk of an \emph{arbitrary} real number $x$ (in the usual sense, with no construction specified for $x$)
in NAFL, for that amounts to quantification, via Proposition~\ref{p3}; however, see also Remark~\ref{arb}.
Nevertheless, indirect quantification over reals may be accomplished by a metamathematical translation of
diagrammatic concepts of Euclidean geometry into operations on sequences of rationals in NPAR. We denote as NRA
the appropriate NAFL metatheory (to be defined below) of NPAR in which this translation is performed.
\begin{definition}\label{d1}
In \emph{NRA}, we associate real numbers with `points' on the real line, which exists a priori as a geometric entity.
Euclidean geometry and its diagrammatic concepts are taken for granted as true. In \emph{NRA}, every mention
of a constant real number~$r$ (e.g.\ $\pi$) \emph{must} be accompanied by a `construction'
for $r$, namely, a specific Cauchy sequence of rationals (that classically converges to $r$); `non-constructive'
existence of real numbers is not permitted in \emph{NAFL} (see Proposition~\ref{p3}).
Any such Cauchy sequence is appropriate in a given instance, and different sequences may be used
in different instances to represent the same real number. In \emph{NRA}, a constant real number $r$
is specified only \emph{after} a Cauchy sequence representing it has been displayed. Define addition, subtraction,
multiplication and division of reals by corresponding termwise operations on the rationals in the Cauchy sequences
representing the reals; we do not permit sequences containing zero terms to represent the denominator of a
quotient of reals. These definitions, adapted from~\cite{sosoa}, are as follows. If $x = \langle q_n \rangle$
and $y = \langle q^{\prime}_n \rangle$ are real numbers, we write $x =_R y$ in \emph{NRA} to mean that
$\lim_n|q_n - q^{\prime}_n| = 0$, i.e., in \emph{NPAR},
\begin{equation*}
\forall \epsilon (\epsilon > 0 \Rightarrow \exists m \forall n (m < n \Rightarrow |q_n - q^{\prime}_n| < \epsilon)),
\end{equation*}
and we write $x <_R y$ in $\emph{NRA}$ to mean that
\begin{equation*}
\exists \epsilon (\epsilon >0 \: \& \: \exists m \forall n (m < n \Rightarrow q_n + \epsilon < q^{\prime}_n))
\end{equation*}
holds in \emph{NPAR}. Similarly, we also make the following translations in \emph{NRA} to formulas of \emph{NPAR}:
\begin{align*}
&x +_R y \to \langle q_n + q^{\prime}_n \rangle \\
&x \cdot_{R} y \to \langle q_n \cdot q^{\prime}_n \rangle \\
&-_Rx \to \langle -q_n \rangle \\
&x /_R y \to \langle q_n / q^{\prime}_n \rangle, \quad  \mbox{where } q^{\prime}_n \ne 0 \\
&0_R \to \langle 0 \rangle; \quad 1_R \to \langle 1 \rangle.
\end{align*}
We have dropped the subscript $Q$ from the relations $\cdot$, $/$, $+$, $-$, $<$ and $>$ on the
rationals, as well as from the symbols $0$ and $1$ (as defined in Appendix~B). We will henceforth drop the
subscript $R$ on these relations and symbols for the reals as well.
\end{definition}
\begin{remark}\label{zbz}
Note that the expression `$\pi / \sqrt{2}$' by itself makes no sense in NAFL unless Cauchy sequences
representing $\pi$ and $\sqrt{2}$ are also specified. The result will be again be a specific Cauchy sequence.
Of course, in this case we know that any such resulting sequence will always represent the same real number.
Consider the expression `$0/0$', where the numerator and denominator are again specified as Cauchy sequences
representing the real number zero. The answer in this case will be dependent on how we specify these sequences,
and may correspond to different real numbers in different instances, or may not converge at all. Whenever the
answer converges, $0/0$ is a well-defined operation in the NAFL version of real analysis. This is so because we
do not require $0/0$ to be a uniquely defined real number like $\pi / \sqrt{2}$; $0/0$
is a (possibly) meaningful real number only \emph{after} Cauchy sequences representing
the zeroes in the numerator and denominator are specified. The definitions of the real zeroes
in the numerator and denominator that \emph{we} provide are only a means for making an axiomatic assertion
of the value of $0/0$, in tune with the Main Postulate.
\end{remark}
We have defined real numbers as Cauchy sequences of rationals representing points on the real line and also the
usual arithmetic operations on these. Next let us consider how one may quantify over specific collections of real
numbers on the real line in NAFL. We will denote such collections as `super-classes', since the reals do not
constitute a class in NAFL (being themselves proper classes).
\begin{definition}\label{d2}
A super-class $I$ of real numbers that has a geometric representation on the real line is defined by specifying
a constant sequence $S$ of rationals that has precisely $I$ as its (super-class of) limit points. Different
sequences $S$ may be used to represent $I$ in different instances. Here $I$ may be a finite collection of reals, an
infinite sequence, or an interval (finite or infinite or semi-infinite), or a combination of these.
\end{definition}
The key point of Definition~\ref{d2} is that only entities $I$ which have a geometric representation on the real
line may be defined as super-classes, as noted below.
\begin{proposition}\label{p4}
Every super-class $I$ of real numbers defined in Definition~\ref{d2} must necessarily be closed, i.e., include
all of its limit points. In particular, open or semi-open intervals of reals, or sequences of reals
excluding their limit points, cannot be represented in \emph{NRA}. Note that such entities do not have a geometric
representation, e.g.\ there is no way to geometrically specify an open/semi-open interval of reals by a diagram.
Similarly, the rationals, when viewed as real numbers in the above geometrical sense, do not constitute a
super-class; any sequence $S$ that includes all rational points represents the entire real line.
\end{proposition}
\begin{proof}
This is immediate from Definition~\ref{d2}. Any sequence $S$ of rationals that has all the elements of an
open/semi-open interval of reals as its limit points \emph{must} also have the end-points of such an interval as
its limit points; this is a classically known result. Similarly, if the sequence $S$ of rationals has an infinite
sequence $\langle \rho_j \rangle$ of reals as its limit points, it must also have the limit points of
$\langle \rho_j \rangle$ as its limit points.
\end{proof}
\begin{remark}
The sequence $S$ of Definition~\ref{d2} has two roles, namely: (a)~As a sequence of
rationals specified by a definite property in NPAR, without any reference to real numbers,
and (b)~As a super-class $I$ of reals when viewed metamathematically from `outside' NPAR, i.e., in NRA. In the latter
case, one must view $S$ as a \emph{simultaneous} specification of \emph{all} of its Cauchy subsequences that represent
the super-class $I$. Thus $S$ satisfies our notion of an infinite class outlined in Remark~\ref{rmcl}. Indeed,
when the theory NPAR \emph{constructs} $S$, it has also unavoidably \emph{constructed} every subsequence of $S$ in
a metamathematical sense. This argument justifies the constructive methods used in the NAFL version of
real analysis~(NRA); there is no other way, given Proposition~\ref{p3}. Cantor's diagonalization argument for the
uncountability of reals is not legal in NRA (and cannot even be formulated in NPAR), because quantification over
reals is not permitted; there is no `cardinality' for a super-class of reals. In fact there is no `list' of reals
in NRA, i.e., there is no way to establish a mapping between any super-class of reals and the natural numbers
$N$ because of the above restriction on quantification.
\end{remark}
\begin{remark}\label{eqr}
Note that Definition~\ref{d2} requires us to extend the meaning of the equality relation $=_R$ in NRA
(see Definition~\ref{d1}) as follows. Two super-classes of reals represented by sequences $\langle q_n \rangle$
and $\langle q^{\prime}_n \rangle$ of rationals are equal if and only if they have precisely the same (super-class
of) limit-points. If $\langle q_n \rangle$ and $\langle q^{\prime}_n \rangle$ happen to be Cauchy sequences, this
modified definition of equality does agree with that in Definition~\ref{d1}. However, the notion of `having the
same super-class of limit points' cannot be formalized within an NAFL theory, such as, NPAR, because it would
require quantification over reals. This is precisely why NRA must remain a metatheory that
appeals to geometric concepts and not a formal NAFL theory of reals (which cannot exist).
\end{remark}
\begin{remark}\label{rat}
The fact that the rationals, when represented metamathematically as reals, do not constitute a super-class in NAFL
(see Proposition~\ref{p4}) means that we cannot consider the rationals as having a separate identity from the
real numbers, under a \emph{geometric} interpretation. The rationals have a purely arithmetical construction in NPAR as
ordered pairs of integers, but they are \emph{not} points on the real line (or solutions of polynomial equations that
permit reals as solutions, etc.) when viewed in this sense. Their construction in NPAR as ordered pairs of integers
is required only for the purpose of defining the reals and such a construction does \emph{not} imply geometric
numberhood for the rationals. When viewed as points on the real line or in any other geometric sense, the rationals
must \emph{always} be viewed as reals, and do not constitute a super-class by themselves. E.g.\ it is not valid in
NAFL to talk of all the `rational points' in the interval [0,1], for no construction is possible for such a
super-class of real numbers. Similarly, integer points on the real line must always be viewed as real numbers.
Any super-class that includes all integer points on the real line must also include $\pm \infty$, where one may
define $\pm \infty$ as real numbers represented by appropriate unbounded sequences of rationals. E.g., an instance
is given below:
\[
\pm \infty =_R \langle \pm_{Q} ((n, 0), (1, 0)) \rangle, \quad \mbox{where } n \in N.
\]
\end{remark}
\begin{remark}\label{arb}
One may wish to represent \emph{arbitrary} real numbers belonging to a super-class in NRA, say, for stating
identities like $x+y =_R y+x$ or $x \nolinebreak =_R \nolinebreak x$. Such a representation is possible if one
provides a \emph{construction} for such an arbitrary real number in the super-class. This is accomplished in
NAFL by using the Main Postulate; an arbitrary real number $x$ in a super-class is taken to be in a superposed
state of assuming all possible values in the super-class until the human mind specifies a
value for it~\cite{1166}. Thus $x$, when unspecified, has the same construction as the super-class itself, i.e.,
as a sequence of rationals with the entire super-class as its only limit points. One may verify that this construction
works when one performs operations like $x + y =_R y + x$, with the usual arithmetical operations
defined similarly to Definition~\ref{d1}. The operations are defined even when we specify $y$ in the above
example while letting $x$ be arbitrary, or let both $x$ and $y$ be arbitrary; however,
note that in these latter instances one needs to interpret $=_R$
as defined in Remark~\ref{eqr}. It must be emphasized that \emph{arbitrary real numbers can only be defined
within legitimate super-classes of NRA}. For example, the following seemingly innocent proposition of classical
real analysis is illegitimate in NRA when $x$ is left arbitrary:
\begin{equation}\label{eqn}
1/x >_R 0 \Rightarrow x >_R 0.
\end{equation}
One might attempt to represent $x$ here as an arbitrary real number in $[-\infty, +\infty]$, which is a legitimate
super-class of NRA. Indeed whenever \emph{we specify} $x$ as a positive or a negative real number (via appropriate
Cauchy sequences of rationals), (\ref{eqn}) does hold. However, when we leave $x$ unspecified, it is in a
superposed state of all possible values as required by NAFL and this includes positive, negative and zero real numbers;
hence $x$ is represented in NRA by the super class $[-\infty, +\infty]$. In this instance, NAFL
would require that both $1/x >_R 0 \; \& \; 1/x \le_{R} 0$ and $x >_R 0 \; \& \; x \le_R 0$ hold in NRA. But then it
follows that $1/x >_R 0 \; \& \; x \le_R 0$ also holds when we leave $x$ unspecified, and this is formally the negation
of (\ref{eqn}). This is the reason for the illegality of (\ref{eqn}) in NRA. Note also that when we specify
$x = 0$ (by an appropriate Cauchy sequence of rationals) in NRA, one could have $1/x = +\infty >_R 0$, which
again results in the falsity of the above assertion. Hence restricting $x$ to the super-class $[0, +\infty]$
of reals also fails to rescue (\ref{eqn}) in NRA. The reason (\ref{eqn}) fails to be legitimate in NRA
is that it attempts to create a super-class $(0, +\infty]$ that does not exist in NRA. However, restricting
$x$ to other super-classes of strictly positive reals (such as, $[1, +\infty]$) will make (\ref{eqn}) legitimate
in NRA, provided $>_R$ is appropriately defined when $x$ is left arbitrary within this super-class. One may also
easily verify that when (\ref{eqn}) is modified to the proposition
\[
1/x \ge_{R} 0 \Rightarrow x \ge_{R} 0,
\]
it can indeed legitimately be asserted as true in NRA when $x$ is restricted to be an arbitrary real number
in $[0, +\infty]$. In this case, one must \emph{define} $\ge_{R}$ such that $x \ge_{R} 0$ is true when
$x$ is left unspecified (and hence is represented by the super-class $[0, +\infty]$).
\end{remark}
\begin{remark}
One may extend the above definitions on the real line to the Euclidean plane or to three-dimensional Euclidean
space by defining points as Cauchy sequences of ordered pairs or triples of rationals that converge to the real
coordinates of the points. One may then define complex numbers, partial derivatives and functions of reals (including
discontinuous ones) that have a geometric representation, by an obvious extension of the results given in this paper.
The arbitrary variables defined in Remark~\ref{arb} are also needed to define functions. For example, to represent
a function as $y = f(x)$, one needs a sequence of pairs of rationals such that it has every pair of values
$(x, f(x))$ as its (super-class of) limit points. Further given any specific $x$ within the super-class representing the
domain, the formula $f(x)$ generates the value $y$ of the function. Finally, if $x$ is left unspecified, its
value is by definition the super-class representing the domain, and the value of $f(x)$ is the super-class
representing the range of the function. The domain and range of $f(x)$ must be legitimate super-classes of NRA.
Discontinuous functions in NRA will take on multiple values at points
of discontinuity, as the classical representation using (semi-) open intervals is not possible; hence these are not
functions in the classical sense. For example, to represent the step function
\[
y = 1, \; \mbox{for} \; x < 0; \quad y = 2, \; \mbox{for} \; x > 0; \quad y = [1, 2] \; \mbox{for} \; x = 0,
\]
one needs a sequence $\langle q_n, q^{\prime}_n \rangle$ of ordered pairs of rationals such that the
(classical) limit points of these sequences are pairs of real numbers of the form $(x, 1)$ for each
$x \le 0$ and $(x, 2)$ for each $x \ge 0$. Further one should also be able to extract the pair $(0,y)$ as limit
points of the above sequence, for each y satisfying $1 \le y \le 2$. The above should, of course, be
the complete super-class of limit points of $\langle q_n, q^{\prime}_n \rangle$. Here we have taken
$y = [1, 2]$ as the super-class representing the `value' of the function at $x = 0$, to conform with the usual
diagrammatic representation of the step function.
\end{remark}
\begin{remark}\label{fp}
Note that \emph{proofs} about operations on real numbers may, in general, be obtained diagrammatically.
Diagrammatic concepts from Euclidean geometry are taken to be consistent and \emph{a priori} true in NRA, before
translation into NPAR. As an example, see~\cite{1255} for a suggested diagrammatic proof of Euclid's fifth
postulate from the first four, in plane geometry. The basic idea of this proof is as follows. Given a line $L$ and
a point $P$ of the plane that is not on $L$, one can deduce from Euclid's first four postulates the existence of
infinitely many line segments of a line $M$, each of which is parallel to $L$; further, this collection of line
segments includes all the points of $M$. In NAFL, this implies the existence of $M$ itself.
In order to translate this proof in NRA, one is faced with a seeming problem: how does one quantify
over line segments, which are themselves super-classes? For example, one needs a representation of infinitely
many line segments of the form $[-1, 1], [-2, 2], [-3, 3], \dots$ in NRA. The appropriate representation
in this case is achieved by two sequences, one of which is a sequence of rationals that has the entire real line
as its limit points and the other is a sequence of pairs of rationals having the limit points
$(-1, 1), (-2, 2), (-3, 3), \dots$ (representing the end-points of the line segments). In particular, the second
sequence must also include $(-\infty, +\infty)$ as its limit point; this is the proof of the fact that
the entire line must necessarily be included in any representation of its line segments that includes all
the points of the line.
\end{remark}
\begin{remark}
In several previous papers~\cite{ifr,1255,ics,acs} we have highlighted the fact that special relativity theory~(SRT)
cannot be formalized as a consistent theory of NAFL. The present paper provides one more such argument, as follows.
SRT requires that a particle can, in principle, acquire \emph{all} velocities strictly less than $c$, the velocity
of light, but the particle velocity can never equal $c$. Formalization of this result requires the existence of the
semi-open interval of reals $[0, c)$. By Proposition~\ref{p4}, such a semi-open interval of reals does not exist
in NRA.
\end{remark}
\begin{definition}\label{der}
The derivative dy/dx in \emph{NRA} may be reduced to a particular means of defining 0/0 as outlined in
Remark~\ref{zbz}. The classical limit of a sequence of reals of the type $\{\Delta y / \Delta x\}$ cannot be
separated from the sequence itself in \emph{NRA} (see Proposition~\ref{p4}); the limit in fact represents $0/0$.
Given $y=f(x)$ and a particular value of $x$, one defines $\Delta x = 0$ as a specific Cauchy sequence of rationals,
and then computes $\Delta y = f(x + \Delta x) - f(x)$ as another Cauchy sequence representing zero. Performing the
division $\Delta y / \Delta x$ yields $dy/dx$ at that value of $x$. Similarly, integration in \emph{NRA}
is an operation of the type $0 \times \infty$.
\end{definition}
\begin{remark}
A detailed exposition of the calculus and the theory of functions of a real variable in NRA will be dealt
with in future work. The derivative, of course, has the geometric interpretation of the slope of a curve.
It is easily seen that essentially every object of Euclidean geometry (whether plane or three-dimensional)
can be translated into NPAR via NRA. As an example of the derivative, consider $y = f(x) = x^2$. Performing
the operation noted in Definition~\ref{der}, one obtains
\begin{align*}
&\Delta y =_R f(x + \Delta x) - f(x) =_R 2x\Delta x + (\Delta x)^2 \\
&\frac{dy}{dx} =_R \frac{\Delta y}{\Delta x} =_R 2x + \Delta x =_R 2x,
\end{align*}
where $\Delta x$ and $x$ are to be replaced by specific Cauchy sequences of rationals representing the real numbers
zero and $x$ respectively. Observe that $x$ may also be left arbitrary within a legitimate super-class of NRA
(see Remark~\ref{arb}). It is extremely important to note that the functional representation $y = f(x)$ must
output a specific Cauchy sequence for $y$ when a Cauchy sequence for $x$ is substituted into $f(x)$. This does
happen in the above example. Note also that the cancellation of $\Delta x$, representing zero, from the numerator
and denominator of $\Delta y / \Delta x$ has no bearing on the final outcome; it is a legal operation in NRA.
\end{remark}
\begin{remark}
The Weierstrass $\epsilon-\delta$ argument for the derivative $dy/dx$ requires the existence
of open intervals of reals and therefore fails in NRA. The NAFL version of real analysis excludes the
paradoxes of classical real analysis, all of which arise from the assumption that non-closed super-classes of reals, in
particular, open/semi-open intervals, exist. As an example, consider one of Zeno's paradoxes of motion, in
which Achilles chases the tortoise and seemingly never catches up with it despite running at a higher velocity.
This paradox only arises when one considers Achilles as reaching, and being confined to, infinitely many
locations (points) at all of which the tortoise must always be ahead. But such a super-class of points does not
exist in NRA; every attempt to represent such a super-class by the NAFL definition
will result in its classical limit point to also be included in it (by Proposition~\ref{p4}), and Achilles does catch
up with the tortoise at this limit point. Similarly the paradoxes of classical measure theory, such as, the
Banach-Tarski paradox, can also be attributed to requiring non-closed super-classes of reals that do not exist in NRA.
\end{remark}
\begin{remark}
It is worth highlighting the fundamental reasons for the paradoxes of classical real analysis and precisely how
they are resolved in NRA. Let us consider again Zeno's paradoxes of motion. The essential paradox here arises
from the fact that infinitely many finite, non-zero intervals of real numbers sum to a finite interval in classical
real analysis. For
example, consider the following infinite sequence of real intervals (in either space or time):
\[
[0, 1/2], \: [1/2, 3/4], \: [3/4, 7/8], \dots
\]
If the above intervals are `stacked' side by side, one gets the real interval $[0, 1]$, reflecting the fact
that $1/2 + 1/4 + 1/8 + \dots = 1$. But note that \emph{each} of the above intervals is finite, non-zero
and non-infinitesimal (infinitesimals do not exist in NRA because NAFL does not permit nonstandard models of NPA).
If there are indeed infinitely many such intervals (as classical real analysis asserts, by mapping these to the
class of natural numbers $N$), then how can their sum be an interval of finite length? For a related paradox,
consider the sequence of nested intervals
\[
[-1, 1], \: [-1/2, 1/2], \: [-1/4, 1/4], \dots
\]
Note that classically, \emph{each} of these nested intervals contains infinitely many (in fact, uncountably many)
points. Then how can the intersection of these intervals contain only a single point, namely, zero? NAFL answers
both of these paradoxes as follows. In each instance, the mapping of the real intervals to $N$ is not legal in
NRA as it amounts to quantification over intervals of reals. When one attempts to construct either of the
above sequences of intervals in NRA by the method noted in Remark~\ref{fp}, one finds that the interval
$[1, 1]$ or $[0, 0]$ respectively of zero length must necessarily get included. This immediately explains
these paradoxes, for it is not possible in NRA to assert that there exist infinitely many finite, non-zero
intervals that sum to a finite interval, as in the above examples. Secondly, it is not even legal in NRA to
ask \emph{how many} intervals exist (including $[1, 1]$ or $[0, 0]$) in each of the above sequences of
intervals, when one constructs them legally in NRA. Such a question assumes that it is legal to quantify over
intervals of reals, whereas it is illegal to even quantify over real numbers in NRA.

Chaitin~\cite{chaitin} cites the paradoxes of classical real analysis and the incompleteness results of G\"odel,
Turing and himself as reasons why the classical real number system must be abandoned. Lynds~\cite{lynds} asserts
that there are no `instants' of time in the real world, which means that Zeno's paradoxes cannot even be formulated.
These authors have essentially rejected the continuum as aphysical; it is interesting to compare their arguments
with the NAFL approach to real analysis presented in this paper.
\end{remark}
\begin{remark}
Exponentiation of reals (or rationals) with rational exponents is relatively straightforward to define.
But exponentiation with real exponents seems problematic in NPAR and will be analyzed in future work.
\end{remark}
\begin{remark}
As an interesting application of real analysis in NAFL, consider Fermat's Last Theorem~(FLT). This may be put in the
form that $r^n \nolinebreak + \nolinebreak s^n \nolinebreak = \nolinebreak 1$ has no solutions in the positive
rationals $r$ and $s$ for integers $n \ge 3$. Let us fix $n$ and $s$ to be specific constants, and consider
a sequence of rationals $\{r_j\}$, such that the limit of this sequence is a real number $\phi$ that solves
$\phi^n \nolinebreak + \nolinebreak s^n \nolinebreak = \nolinebreak 1$; FLT requires that no member of the
sequence $\{r_j\}$ should solve this equation. But we have seen from Remark~\ref{rat} that the rationals,
when substituted into polynomial equations that permit reals as solutions, \emph{must} be considered as reals.
So FLT requires that a super-class of reals $\{r_j\}$ must exist that fails to solve the above equation, while
its limit point does solve it. This amounts to requiring that the limit point of this purported super-class must be
excluded from it, but this is not possible by Proposition~\ref{p4}. One concludes that the truth of FLT cannot be
meaningfully represented if the geometric interpretation of real numbers is imposed upon Peano Arithmetic~(NPA).
Therefore NAFL requires that any proof (or refutation) of FLT must be carried out entirely within NPA, without
invoking real numbers. This reinforces our earlier conclusion of Remark~\ref{rmnon} that consistency of NPA
demands its completeness.
\end{remark}
\begin{remark}
Simpson~\cite{sosoa} (with his subsystems of second-order arithmetic), Feferman~\cite{fefer,fefer2} (predicativism) and
Weaver~\cite{nik,nik2} (mathematical conceptualism) have all attempted to restrict classical infinitary reasoning
so as to eliminate the set-theoretic paradoxes. The fundamental difference between
the NAFL approach of this paper and those of the above authors is that NAFL does not accept the existence of infinite
sets and quantification over (infinite) proper classes. Neither does NAFL accept the notion of an arbitrary infinite
class (as a free variable), which amounts to quantification in NAFL. The point we wish to make here is that
quantification over proper classes is fundamentally and unavoidably part of infinitary reasoning and so stands
rejected by NAFL, as noted in the proof of Proposition~\ref{p3}. From the NAFL point of view,
we dispute Simpson's~\cite{sosoa} claim that he has achieved partial realizations of Hilbert's program (to justify
classical infinitary reasoning from the finitary standpoint) because Simpson has really used infinitary methods.
In fact NAFL shows that Hilbert's program is decisively settled negatively~(see also \cite{1166}); classical
infinitary reasoning stands refuted from the strictly finitary standpoint developed in NAFL. We also mention in
passing that Weaver~\cite{nik2} has severely criticized Feferman's~\cite{fefer2} approach to predicativism.
We honestly believe though, that neither Weaver's nor Feferman's approach is predicativism in the strict sense; they
accept the existence of infinite sets of natural numbers, which are essentially impredicative objects from the
NAFL standpoint~(see Sect.~\ref{infsets}). In fairness to these authors, they have characterized their approach
as predicativism `given the natural numbers'.
\end{remark}

\section{Foundations of Computability Theory in NAFL}
Here we will outline the basic arguments in \cite[Sects.~6~and~7]{acs}, and establish the connection with real analysis
in NAFL. A Turing machine~(TM), \emph{by definition}, must either halt or not halt. If $P$ is the proposition that
a given TM halts, then $P \vee \neg P$ is unavoidably built into the definition of that TM. It follows
that $P$ cannot be undecidable in any consistent NAFL theory~T in which the existence of that TM is formalized;
the non-classical model for T required by Proposition~\ref{p1} in which $P \& \neg P$ is the case cannot exist,
for T must prove $P \vee \neg P$ (and hence, either $P$ or $\neg P$). Any infinite (proper) class in an NAFL theory
must be recursive; whether an object belongs to or does not belong to that class cannot be undecidable in a
consistent NAFL theory because such undecidability will violate the uniqueness required by the axiom for
extensionality for classes, which is an essential ingredient of NAFL theories in which infinite classes
exist~(e.g.\ see Sect.~\ref{infsets}, as well as Remarks~\ref{rmcl}~and~\ref{ae}). Note that the above arguments
require the non-classical features of NAFL, namely, Proposition~\ref{p1} and the Main Postulate. It is clear
that a new paradigm for computability theory is required in NAFL.
\begin{proposition}[NAFL Computability Thesis, NCT]\label{p5}
Every infinite (recursive, proper) class that exists in \emph{NAFL} theories must be effectively computable,
i.e., there exists an algorithm that computes it. Here `recursive' or `computable' does not necessarily mean
`computable by a classical Turing machine'.
\end{proposition}
\begin{proof}
It is clear that for every element of an infinite class, which is a finite set, there exists an algorithm that
computes it, i.e., it is effectively computable. Note that nonstandard models of arithmetic do not exist in
NAFL; see the proof of Proposition~\ref{p2} and Remark~\ref{rmnon}. Therefore the existence of infinitely many
algorithms that compute every (standard finite) element of an infinite class implies that the entire infinite class
has been computed by these algorithms. But by Remark~\ref{rmcl}, these algorithms cannot operate one or finitely
many at a time and complete the computation process; by induction, there will always exist infinitely many
elements remaining to be computed. One concludes that in NAFL, the existence of infinitely many algorithms
as noted above necessarily implies the existence of \emph{an} algorithm~$\mathcal{A}$ that computes the infinite
class. One may think of $\mathcal{A}$ as operating on the infinitely many algorithms corresponding to all elements
and `parallelizing' them, so that the infinite class gets computed \emph{simultaneously}, rather than
finitely many at a time, as required by the NAFL interpretation of an infinite class (see Remark~\ref{rmcl}). Note
that the nonexistence of $\mathcal{A}$ would leave us with a paradox in NAFL, namely, that a class of
algorithms collectively managed to compute an infinite class, although each of them did only a finite computation.
The above proof of the existence of $\mathcal{A}$ required the non-classical concepts of NAFL,
namely, the Main Postulate, Proposition~\ref{p1} and the nonexistence of nonstandard integers.
The proof that every infinite class must be recursive also required these concepts. Hence $\mathcal{A}$
need not be a classical algorithm, i.e., the infinite class in question need not be computable by a classical
Turing machine.
\end{proof}
\begin{remark}\label{qalg}
It is a theorem of classical recursion theory that every recursive class must be computable by a Turing machine.
But many of these standard results fail in NAFL, whose concepts of `computability'
and `algorithm' must necessarily be different from the classical concepts.
Proposition~\ref{p5} implies that Turing's halting routine $H$ must exist in NAFL.
We believe that such an algorithm must be non-classical, not because of Turing's argument (which fails in NAFL),
but because requiring $H$ to be classical would make it a self-referential entity that would lead to contradictions
in NAFL. Thus NAFL permits hypercomputation, and there is a case for believing that $H$ must be a quantum algorithm, for
the following reasons. NAFL also justifies quantum superposition and entanglement~\cite{ijqi,1166,acs}, which are
non-classical and metamathematical phenomena, i.e., they are confined to the metatheory (semantics)
of NAFL theories. Hence a purported quantum algorithm for $H$, which only computes the halting
decisions for classical Turing machines, would reside in the metatheory and will not be a self-referential entity.
Secondly, there is evidence that quantum algorithms permit hypercomputation~\cite{kieu} and infinite
parallelism~\cite{ziegler}. The latter feature may also be justified by noting that in NAFL, a quantum algorithm
is permitted to access and compute a truly random element of an infinite class. When such an element
is not specified, it must be in a superposed state of assuming all values in the infinite class~(see Remark~\ref{arb})
and this corresponds to infinite parallelism.
\end{remark}
\begin{remark}
An interesting analogy is suggested between the rational/real numbers on the one hand,
and classical/quantum algorithms on the other. It is known that for every classical algorithm, there
is a quantum algorithm that achieves the same result; the converse of this assertion is controversial
and not necessarily true. The classical and quantum representations of an algorithm (when they both exist)
may be thought of as corresponding to the dual representation of a number as a rational and real respectively.
As noted in Remark~\ref{rat}, the rational representation does not have the geometrical significance
of a point on the real line; in a similar sense, the classical representation of an algorithm merely
encodes a meaningless finite string of symbols, rather than an algorithm that executes. When thought
of as an algorithm, the quantum representation \emph{must} always be used. This is justified
by the fact that the infinitely many classical algorithms (in the execution mode) by themselves do
not constitute a class because they cannot be separated from the quantum algorithms that `parallelize'
them (see Proposition~\ref{p5}, its proof and Remark~\ref{qalg}). Quantum algorithms may be
thought of as infinite sequences of the classical representations; these sequences
are only needed to define the quantum algorithms (in the same sense that reals are infinite sequences
of rational representations that are only needed to define the reals). Such sequences of classical
representations may also be used to define `super-classes' of quantum algorithms, in the same manner that NRA
defines super-classes of reals. An executing algorithm, even when it halts and has a classical representation,
is to be thought of as an infinite object corresponding to its quantum representation. This may be
because the machine that executes the algorithm has an infinite tape, for example; the instruction `halt'
may be interpreted as a tacit requirement not to operate further on this infinite tape. Secondly, the classical model of
computation is not valid in NAFL and the algorithm may at some point access infinitely many execution states at the
same time. A quantum algorithm that does not have a classical representation may be thought of as corresponding to an
irrational real number. In conclusion, we observe that important potential applications of the NAFL paradigms
for real analysis and computability theory are in the areas of quantum and autonomic computing~\cite{acs}.
\end{remark}

\section*{Appendix~A. Proof of Proposition~\ref{p1}}
\begin{proof}
By the Main Postulate of NAFL, $P$~($\neg P$) can be the case in T if and only if
$P$~($\neg P$) has been asserted \emph{axiomatically}, by virtue of its
provability in T*. In the absence of any such axiomatic assertions (\emph{e.g.}~if T*=T),
it follows that neither $P$ nor $\neg P$ can be the case in T and hence
$P \vee \neg P$ cannot be a theorem of T. The classical refutation of
$P \& \neg P$ in T proceeds as follows: `If $P$~($\neg P$) is the case, then
$\neg P$~($P$) cannot be the case', or equivalently, `$\neg P$~($P$) contradicts
$P$~($\neg P$)'. But, by the Main Postulate, this argument fails in NAFL and
amounts to a refutation of $P \& \neg P$ in T*=T+$P$~(T+$\neg P$), and \emph{not}
in T as required. Careful thought will show that the classical refutation of
$P \& \neg P$ in T is the \emph{only possible} reason for $\neg (P \& \neg P)$
to be a theorem of T, and it fails in NAFL. The intuitionistic refutation of
$P \& \neg P$ in T is flawed and also fails in NAFL, as will be shown in the ensuing paragraph. By the
completeness theorem of FOPL (which, as noted in Remark~\ref{rmnon}, is taken for granted as a metamathematical
principle in NAFL) it follows that there must exist a non-classical model for T in which $P \& \neg P$ is satisfiable.
\end{proof}

Consider the law of non-contradiction as stated in a standard system of intuitionistic
first-order predicate logic due to S.~C.~Kleene, namely, \newline
$\neg P \Rightarrow (P \Rightarrow Q)$. This formula asserts that from contradictory premises
$P$ and $\neg P$, an \emph{arbitrary} proposition $Q$ can be deduced,
which is absurd. Hence $\neg (P \& \neg P)$ seemingly follows. However, note that in intuitionism, truth is
provability (not necessarily in a specific theory T); together with
the intuitionistic concept of negation, it follows that an assertion of
$\neg (P \& \neg P)$ is the same as deducing an absurdity from $P \& \neg P$, or
equivalently, from contradictory premises $P$ and $\neg P$. But we have seen
that the `absurdity' referred to here is precisely the fact
that \emph{any proposition can be deduced}, given contradictory premises! The above
`proof' of $\neg (P \& \neg P)$ from contradictory premises,
mandated by the intuitionistic concepts of truth and negation, is flawed
because \emph{any proposition can be so deduced}. Note that this `proof' is formally
indistinguishable from one in which $\neg (P \& \neg P)$ is \emph{substituted} for the
deduced arbitrary proposition~$Q$. In NAFL, it is not possible to deduce an arbitrary proposition
from contradictory premises~\cite{635} in a non-classical model, and so the flawed
intuitionistic argument for $\neg (P \& \neg P)$ fails in any case. Indeed, as explained
in \cite{635}, the argument for deducing an arbitrary proposition would normally
proceed as follows:
\begin{eqnarray}
&& P \& \neg P \Rightarrow P, \nonumber  \\
&&P \Rightarrow P \vee Q,  \nonumber    \\
&&P \& \neg P \Rightarrow \neg P,  \nonumber  \\
&&\neg P \&  (P \vee Q) \Rightarrow Q. \nonumber
\end{eqnarray}
The final step fails in a non-classical NAFL model for a theory~T (in which $P \& \neg P$
is the case) because this step presumes $\neg (P \& \neg P)$.

\section*{Appendix~B. The formal system NPA and its extension NPAR}
We will outline the formal system NPA (i.e., the NAFL version of first-order Peano Arithmetic) along the lines
of the description given in Chapter~1 of \cite{sosoa}, suitably modified for our purposes. Throughout, a
natural deduction system of classical first order predicate logic with equality~(FOPL) is assumed. The language of
NPA admits number variables, denoted by $i, j, k, m, n, \dots,$ which are intended to range over the class $N$ of
all natural numbers $\{0, 1, 2, \dots \}$. Capital letters $L, M, N, \dots,$ are reserved for constant
(infinite) classes; there are no class variables, and we do not need finite classes for our purposes (though these
can be added if required). These infinite classes are defined by the notation $\{n: \phi(n)\}$, for each
L-formula $\phi(n)$ that holds for infinitely many values of $n$.

The terms and formulas of the language~(L) of NPA are as follows. Numerical terms are number variables, the constant
symbols $0$ and $1$, and $t_1 + t_2$ and $t_1 \cdot t_2$ whenever $t_1$ and $t_2$ are numerical terms. Here $+$ and $.$
are binary operation symbols intended to denote addition and multiplication of natural numbers. Numerical terms
are intended to denote natural numbers. Atomic formulas are $t_1 = t_2$, $t_1 < t_2$ and $t_1 \in X$ where $t_1$
and $t_2$ are numerical terms and $X$ is any (constant) infinite class. The intended meanings of these atomic
formulas are that $t_1$ equals $t_2$, $t_1$ is less than $t_2$, and $t_1$ is an element of $X$. Formulas are
built up from atomic formulas by means of propositional connectives $\&$, $\vee$, $\neg$, $\Rightarrow$,
$\Leftrightarrow$ (and, or, not, implies, if and only if) and universal (number) quantifiers $\forall n$,
$\exists n$ (for all $n$, there exists $n$). There are no class quantifiers. A sentence is a formula with no
free variables. In writing terms and formulas of L, parentheses and brackets will be used to indicate grouping
and some obvious abbreviations including the symbols $\notin$ and $\ne$ and the numbers $2, 3, \dots$ will also
be used for convenience.

The axioms of NPA are the following L-formulas (universal quantification of number variables is assumed):
\begin{align*}
\intertext{(i) basic axioms}
&n+1 \ne 0  \\
&m+1 = n+1 \Rightarrow m=n  \\
&m+0 = m  \\
&m+(n+1) = (m+n)+1  \\
&m \cdot 0 = 0  \\
&m \cdot (n+1) = (m \cdot n) + m  \\
&\neg \: m < 0  \\
&m < n+1 \Leftrightarrow (m < n \vee m=n)
\intertext{(ii) induction axiom scheme (for the L-formulas $\phi(n)$)}
&(\phi(0) \; \& \; \forall n (\phi(n) \Rightarrow \phi(n+1))) \Rightarrow \forall n \phi(n)
\end{align*}
Note that the L-formulas $\phi(n)$ may contain independently defined constant class symbols. Further,
by Remark~\ref{rmnon}, $\phi(n)$ must be a sentence; the free (number) variables are automatically assumed to be
universally quantified in NAFL. By Proposition~\ref{p2}, NPA must prove the existence of every infinite
class $\mathcal{C}$ generated by L-formulas, as indicated by the following theorem scheme, which we
will denote as the theorem scheme of comprehension:
\[
\forall n (n \in \mathcal{C} \Leftrightarrow \phi(n)).
\]
Here $\mathcal{C}$ is the generic notation that stands for $\{n: \phi(n)\}$. The theorem scheme applies for
each L-formula $\phi(n)$ that holds for infinitely many values of $n$. Note that $\phi(n)$ may not contain
$\mathcal{C}$, but could possibly contain other (constant) class symbols that have been defined independently
of $\mathcal{C}$. The given instance of the theorem scheme says that there exists an infinite class
$\mathcal{C} = \{n: \phi(n)\}$, which is the class of all $n$ such that $\phi(n)$ holds. As noted in
Remarks~\ref{rmcl}~and~\ref{ae}, the `infinite class' referred to here is always identified by all of its elements,
i.e., it is always a `class' in the extensional sense and NPA must also prove the axiom of extensionality
for classes in those instances where infinite classes are involved.

The above completes the description of the language~L and the axioms of NPA, whose objects are the natural numbers.
The theorems of NPA are deduced from the axioms using the classical rules of inference. The restrictions required
by NAFL as noted in Sect.~\ref{intro}, and in particular, Sects.~\ref{mp}~and~\ref{tsps}, apply to NPA.

Next consider the extension NPAR of NPA in order to handle integers, rationals and sequences of rationals.
The language LR of NPAR augments L, firstly by admitting ordered pairs of natural numbers in the
form $(m, n)$, which are the integers, belonging to the infinite class $Z$. The axioms for the integers are
as follows (universal quantification over the number variables is assumed):

\begin{align*}
\intertext{(i) equality}
&(m, n) =_Z (i, j) \Leftrightarrow m+j = n+i  \\
\intertext{(ii) addition}
&(m, n) +_Z (i, j) =_Z (m+i, n+j)  \\
\intertext{(iii) negative integers}
&-_Z(m, n) =_Z (n, m)  \\
\intertext{(iv) multiplication}
&(m, n) \cdot_{Z} (i, j) =_Z (mi+nj, mj+ni)  \\
\intertext{(v) zero and one}
&0_Z =_Z (0, 0); \quad  1_Z =_Z (1, 0)  \\
\intertext{(vi) order}
&(m, n) <_Z (i, j) \Leftrightarrow m + j < n + i
\end{align*}

Note that we have avoided mentioning equivalence classes while defining the integers. Although
equivalence classes of pairs of naturals defining specific integers do exist in NPAR (via the theorem of
comprehension), one cannot quantify over these in NAFL because they are infinite classes. Hence it is
not possible to use equivalence classes to define infinitely many integers in NPAR.
Following Simpson~\cite[Chap.~1]{sosoa}, one may also define $Z$ to be a class of \emph{representatives}
of the appropriate equivalence classes, e.g.:
\[
Z=\{(0,0), (1,0), (0,1), (2,0), (0,2), (3,0), (0,3), \dots \}.
\]
But note that \emph{formally}, the above definition of $Z$ does not require any mention of equivalence classes.
Thus in NPAR, whenever we formulate a proposition to hold for all integers, it is to be understood that the
proposition is true of all members of $Z$ in the above form, or alternatively, for all ordered pairs
$(m, n)$ as noted in the above axioms.

The rationals (denoted by the infinite class $Q$) are defined in NPAR as ordered pairs $(a, b)$ of integers,
where $b$ is restricted to be a positive integer. The axioms for rationals are as follows (universal
quantification over free variables is assumed).

\begin{align*}
\intertext{(i) equality}
&(a, b) =_Q (c, d) \Leftrightarrow a \cdot_{Z} d =_Z b \cdot_{Z} c  \\
\intertext{(ii) addition}
&(a, b) +_Q (c, d) =_Q (a \cdot_{Z} d +_Z b \cdot_{Z} c, b \cdot_{Z} d)  \\
\intertext{(iii) negative rationals}
&-_Q(a, b) =_Q (-_{Z}a, b)  \\
\intertext{(iv) multiplication}
&(a, b) \cdot_{Q} (c, d) =_Q (a \cdot_{Z} c, b \cdot_{Z} d)  \\
\intertext{(v) division}
&(a, b) /_{Q} (c, d) =_Q (a \cdot_{Z} d, b \cdot_{Z} c), \mbox{where } c \ne 0_Z
\intertext{(vi) zero and one}
&0_Q =_Q (0_Z, 1_Z); \quad  1_Q =_Q (1_Z, 1_Z)  \\
\intertext{(vii) order}
&(a, b) <_Q (c, d) \Leftrightarrow a \cdot_{Z} d <_Z b \cdot_{Z} c
\end{align*}

Note that the axiom for division is redundant, but it will be useful for our purposes in defining division
of real numbers. Following Simpson~\cite[Chap~1]{sosoa}, one may also define $Q$ in NPAR as a class of
representatives of equivalence classes; these representatives may be chosen to be minimal in the sense that all
common factors are cancelled out. Again, $Q$ must be constructively defined in NPAR without any reference to
equivalence classes in order to avoid quantification over these.

The idea behind the above definitions of $Z$ and $Q$ within NPAR is that $(m, n)$ corresponds to
the integer $m - n$, while $(a, b)$ (with $b$ restricted to be a positive integer) corresponds to the
rational $a/b$.

A sequence of rational numbers is defined to be a function $f : N \to Q$. We denote such a sequence
as $\langle q_n: n \in N\rangle$ or simply $\langle q_n \rangle$, where $q_n = f(n)$. A Cauchy sequence of
rational numbers is a sequence $\langle q_n: n \in N \rangle$ such that
\[
\forall \epsilon(\epsilon > 0 \Rightarrow \exists m \forall n (m < n \Rightarrow |q_m - q_n| < \epsilon)).
\]
Here $\epsilon$ ranges over $Q$ and we have dropped the obvious subscripts on the operations $>$, $<$,
$-$, etc. Propositions~\ref{p2}~and~\ref{p3} require that NPAR must prove the existence of every sequence of
rationals (including Cauchy sequences) definable in the language of NPAR. Remarks~\ref{rmcl}~and~\ref{ae}
apply. Of course, the language of NPAR must also admit sequences of integers definable by its formulas,
and NPAR must prove the existence of these as well.

\setlength{\unitlength}{1cm}
\begin{picture}(10,0.1)
\put(-0.5,0){\line(1,0){5}}
\end{picture}
\small \\
IBM and Autonomic Computing are trademarks of the International Business Machines Corporation
in the United States, other countries, or both.


\begin{thebibliography}{99}
\bibitem{ijqi}
R.~Srinivasan, The quantum superposition principle justified in a new non-Aristotelian finitary logic,
\textit{International Journal of Quantum Information}, Vol.~3, No.~1~(2005), pp. 263-267. Presented at the international
conference \textit{Foundations of Quantum Information} (FQI04), University of Camerino, Italy, April 16-19, 2004.
Eprint available at: http://philsci-archive.pitt.edu/archive/00001923/ .

\bibitem{1166}
R.~Srinivasan, Platonism in classical logic versus formalism in the proposed non-Aristotelian finitary logic,
Eprint, May~2003, available at: http://philsci-archive.pitt.edu/archive/00001166/ .

\bibitem{bcp}
R.~Srinivasan, Logical analysis of the Bohr Complementarity Principle in Afshar's experiment under the
NAFL interpretation, arXiv: quant-ph/0504115 .

\bibitem{acs}
R.~Srinivasan and H.~P.~Raghunandan, On the existence of truly autonomic computing systems and the link
with quantum computing, arXiv: cs.LO/0411094 .

\bibitem{635}
R.~Srinivasan, Quantum superposition justified in a new non-Aristotelian finitary logic,
Eprint, May~2002, available at: http://philsci-archive.pitt.edu/archive/00000635/ .

\bibitem{Mcr}
M.~Machover, \emph{Set theory, logic and their limitations}, Cambridge University Press, Cambridge, 1996.

\bibitem{Coh}
P.~J.~Cohen, \emph{Set theory and the continuum hypothesis}, W.~A.~Benjamin, Inc., New York, 1966.

\bibitem{sosoa}
S.~G.~Simpson, \emph{Subsystems of Second Order Arithmetic}, Second Edition, 2005. To be published by the Association
For Symbolic Logic in their book series \emph{Perspectives in Logic}. Chapter one of this book is available at:
http://www.math.psu.edu/simpson/sosoa/ .

\bibitem{ifr}
R.~Srinivasan, Inertial frames, special relativity and consistency,
Eprint, June~2002. http://philsci-archive.pitt.edu/archive/00000666/ .

\bibitem{ics}
R.~Srinivasan, Relativistic determinism: the clash with logic, Poster presentation at the
\emph{International Conference on the Ontology of Spacetime}, Concordia University,
Montreal, May~11--14 (2004). Available at: \\
http://alcor.concordia.ca/$\sim$scol/seminars/conference/abstracts.html .

\bibitem{fefer}
S.~Feferman, Why a little bit goes a long way: Logical foundations of scientifically applicable mathematics,
in \emph{PSA~1992}, Vol.~II, pp.~442-455, 1993. Reprinted as Chapter~14 in \emph{In the light of logic: Logic
and Computation in Philosophy}, Oxford University Press, 1998.

\bibitem{fefer2}
S.~Feferman, Systems of predicative analysis, \emph{J. Symbolic Logic}, Vol.~29 (1964), pp.~1-30.

\bibitem{nik}
N.~Weaver, Mathematical conceptualism, arXiv: math.LO/0509246.

\bibitem{nik2}
N.~Weaver, Predicativity beyond $\Gamma_0$, arXiv: math.LO/0509244.

\bibitem{1255}
R.~Srinivasan, On the logical consistency of special relativity theory and non-Euclidean geometries: Platonism
versus formalism, July~2003, available at: http://philsci-archive.pitt.edu/archive/00001255/ .

\bibitem{chaitin}
G.~J.~Chaitin, How real are real numbers?, arXiv: math.HO/0411418.

\bibitem{lynds}
P.~Lynds, Time and classical and quantum mechanics: indeterminacy versus discontinuity.
\emph{Foundations of Physics Letters}, Vol.~16, Number 4 (2003), pp.~343-355.

\bibitem{kieu}
T.~Kieu, An anatomy of a quantum adiabatic algorithm that transcends the Turing computability,
\textit{International Journal of Quantum Information}, Vol.~3, No.~1~(2005), pp.~177-182.

\bibitem{ziegler}
M.~Ziegler, Does quantum mechanics allow for infinite parallelism? arXiv: quant-ph/0410141 .

\end{thebibliography}
\end{document}